\documentclass[11pt]{article}
\usepackage{amssymb}
\usepackage{stmaryrd}
\usepackage{amssymb}
\usepackage{amssymb}
\usepackage{amssymb}
\usepackage{amssymb}
\usepackage{amssymb}
\usepackage{amssymb}
\usepackage{amssymb}
\usepackage{amssymb}
\usepackage{amssymb}
\usepackage{amssymb}
\usepackage{amssymb}
\usepackage{amssymb}
\usepackage{amssymb}
\usepackage{amssymb}
\usepackage{amssymb}
\usepackage{amssymb}
\usepackage{amssymb}
\usepackage{amssymb}
\usepackage{amssymb}
\usepackage{amssymb}
\usepackage{amssymb}
\usepackage{amssymb}
\usepackage{amssymb}
\usepackage{amssymb}
\usepackage{amssymb}
\usepackage{amssymb}
\usepackage{amssymb}
\usepackage{amssymb}
\usepackage{amssymb}
\usepackage{amssymb}
\usepackage{amssymb}
\usepackage{amssymb}
\usepackage{amssymb}
\usepackage{amssymb}
\usepackage{amssymb}
\usepackage{amssymb}
\usepackage{amssymb}
\usepackage{indentfirst}
\usepackage{amscd}
\usepackage{amsmath}
\usepackage{latexsym}
\usepackage{amsfonts}
\usepackage{amssymb}
\usepackage{amsthm}
\usepackage{mathrsfs}
\usepackage{ams,amsmath,amsthm,amsfonts,amssymb,bm}

\begin{document}
\title{Further development of positive semidefinite solutions of the
operator equation $\sum_{j=1}^{n}A^{n-j}XA^{j-1}=B$
 \thanks
{ This work is supported by National Natural Science Fund of China
(10771011 and 11171013
).}}
\author{ Jian Shi$^{1}$ \ \ Zongsheng Gao \hskip 1cm  \\
{\small   LMIB $\&$ School of Mathematics and Systems Science,}\\
{\small Beihang University, Beijing, 100191, China}\\ }
         \date{}
         \maketitle

\maketitle \baselineskip 16pt \line(1,0){420}

\noindent{\bf Abstract}\ \ \ In \cite{Positive semidefinite
solutions}, T. Furuta discusses the existence of positive
semidefinite solutions of the operator equation
$\sum_{j=1}^{n}A^{n-j}XA^{j-1}=B$. In this paper, we shall apply
Grand Furuta inequality to study the operator equation.
A generalized special type of $B$ is obtained
 due to \cite{Positive semidefinite solutions}. \vspace{0.2cm}

\noindent {\bf Keywords}:  Furuta inequality; Grand Furuta
inequality; operator equation; matrix equation; positive
semidefinite operator; positive definite operator

\noindent {\bf Mathematics Subject Classification:} 15A24, 47A62,
47A63. \vspace{0.2cm} \vskip   0.2cm \footnotetext[1]{Corresponding author. 

E-mail
addresses: shijian@ss.buaa.edu.cn\ (J. Shi), zshgao@buaa.edu.cn\ (Z. Gao).}

\setlength{\baselineskip}{20pt}
\section{Introduction }
\indent A capital letter $T$ means a bounded linear operator on a
Hilbert space. $T \geqslant 0$ and $T>0$ mean a positive
semidefinite operator and a positive definite operator,
respectively.

In the middle of last century, E. Heinz et al. studied
operator theory and obtained the following famous theorem:

\noindent {\bf Theorem LH (L\"{o}wner-Heinz inequality,
\cite{Lowner}  \cite{Heinz}).} If $A\geqslant B\geqslant 0$, then
$A^{\alpha}\geqslant B^{\alpha}$ for any $\alpha\in [0, 1]$.

\indent In 1987, T. Furuta proved the following operator inequality
as an extension of Theorem LH:

\noindent {\bf Theorem F (Furuta inequality, \cite{Furuta
inequality}).} If $A\geqslant B\geqslant 0$, then for each
$r\geqslant 0$,
$$
(B^{\frac r 2}A^{p}B^{\frac r 2})^{\frac 1 q}\geqslant (B^{\frac r
2}B^{p}B^{\frac r 2})^{\frac 1 q} \eqno (1.1)
$$
$$
(A^{\frac r 2}A^{p}A^{\frac r 2})^{\frac 1 q}\geqslant (A^{\frac r
2}B^{p}A^{\frac r 2})^{\frac 1 q} \eqno (1.2)
$$
hold for $p\geqslant 0$, $q\geqslant 1$ with $(1+r)q\geqslant p+r$.

\indent K. Tanahashi, in \cite{Best Furuta inequality}, proved the
conditions $p$, $q$ in Theorem F are best possible if $r\geqslant
0$.

In 1995, T. Furuta showed another operator inequality which
interpolates Theorem F:

\noindent {\bf Theorem GF (Grand Furuta inequality, \cite{GF
inequality}).} If $A\geqslant B\geqslant 0$ with $A>0$, then for
each $t\in [0, 1]$ and $p\geqslant 1$,
$$
A^{1-t+r}\geqslant \{A^{\frac r 2}(A^{-{\frac t 2}}B^{p}A^{-{\frac t
2}})^{s}A^{\frac r 2}\}^{\frac {1-t+r}{(p-t)s+r}} \eqno (1.3)
$$
holds for $s\geqslant 1$ and $r\geqslant t$.

Afterwards, some nice proof of Grand Furuta inequality are shown,
such as \cite{Mean theoretic}, \cite{Simplified proof GF}. K.
Tanahashi, in \cite{Best GF inequality}, proved that the outer
exponent value of (1.3) is the best possible. Later on, the proof
was improved, see \cite{A short proof best GF}, \cite{Simplified
proof Best GF}.

Recently, T. Furuta proved the following theorem by Furuta
inequality:

\noindent {\bf Theorem A (\cite{Positive semidefinite solutions}).}
If $A$ is a positive definite operator and $B$ is positive
semidefinite operator. Let $m$ and $n$ be nature numbers. There
exists positive semidefinite operator solution $X$ of the following
operator equation:
$$
\sum^{n}_{j=1}A^{n-j}XA^{j-1}=A^{\frac
{nr}{2{(m+r)}}}(\sum^{m}_{i=1}A^{\frac {n(m-i)}{m+r}}BA^{\frac
{n(i-1)}{m+r}})A^{\frac {nr}{2(m+r)}} \eqno (1.4)
$$
for $r$ such that $ \begin{cases}
        r\geqslant 0 ,  & \text{if }  n\geqslant m ;\\
        r\geqslant {\frac {m-n}{n-1}} , & \text{if }  m\geqslant n\geqslant 2 .
\end{cases}$
\\

In the rest of this short paper, we
shall apply Grand Furuta inequality to
 discuss the existence of positive semidefinite solution of operator equation
$\sum^{n}_{j=1}A^{n-j}BA^{j-1}=B$, and show
a generalized special type of $B$ due to Theorem A.

\section{Extension of Furuta's result }

\noindent{\bf Lemma 2.1 (\cite{Bhatia}, \cite{Positive semidefinite solutions}).}
Let $A$ be a positive definite operator and $B$ a positive
semidefinite operator. Let $m$ be a positive integer and $x\geqslant
0$. Then $ \left. {\frac
{d}{dx}}[(A+xB)^{m}] \right| _{x=0}
=\sum_{j=1}^{m}A^{m-j}BA^{j-1}$.

\noindent{\bf Theorem 2.1.} Let $A$ be a positive definite operator
and $B$ be a positive semidefinite operator. Let $m$, $n$, $k$ be
positive integers, $t\in[0, 1]$. There exists positive semidefinite
operator solution $X$ which satisfies the operator equation:
\begin{equation}
\begin{split}
\sum^{n}_{j=1}A^{n-j}XA^{j-1}=& A^{\frac
{nr}{2(m-t)k+2r}}\{\sum^{k}_{i=1}A^{\frac
{n(m-t)(k-i)}{(m-t)k+r}}[A^{{-{\frac t 2}}\cdot{\frac
{n}{(m-t)k+r}}}\cdot(\sum^{m}_{j=1}A^{\frac
{n(m-j)}{(m-t)k+r}}BA^{\frac
{n(j-1)}{(m-t)k+r}})\cdot \\
&A^{{-{\frac t 2}}\cdot{\frac {n}{(m-t)k+r}}}]A^{\frac
{n(m-t)(i-1)}{(m-t)k+r}}\}A^{\frac {nr}{2(m-t)k+2r}}
\end{split}
\end{equation}
for $r$ such that $ \begin{cases}
        r\geqslant t ,  & \text{if }\   (1-t)n\geqslant (m-t)k \  ;\\
        r\geqslant max \{{\frac {(m-t)k-(1-t)n}{n-1}} , t\} , & \text{if }
        \   (m-t)k\geqslant (1-t)n \  \text{with } n\geqslant 2 \ .
\end{cases}$
\\

\noindent{\bf Proof.} As in the proof of
 [\cite{Positive semidefinite solutions}, Theorem 2.1], by
$A+xB\geqslant A>0$ holds for any $x\geqslant 0$, then
$A^{-1}\geqslant (A+xB)^{-1}>0$. Replace $A$ by $A^{-1}$, $B$ by
$(A+xB)^{-1}$, $p$ by $m$, $s$ by $k$ in (1.3), and take reverse, we
have
$$
(A^{\frac r 2}(A^{{-}{\frac t 2}}(A+xB)^{m}A^{{-}{\frac t
2}})^{k}A^{\frac r 2})^{\frac {1-t+r}{(m-t)k+r}}\geqslant A^{1-t+r}.
\eqno (2.2)
$$

For any $\alpha \in [0, 1]$, apply L\"{o}wner-Heinz inequality to
(2.2), and take ${\frac 1 n}={\frac {1-t+r}{(m-t)k+r}}\cdot \alpha$,
the following inequality is obtained:
$$
(A^{\frac r 2}(A^{-{\frac t 2}}(A+xB)^{m}A^{-{\frac t
2}})^{k}A^{\frac r 2})^{\frac 1 n}\geqslant A^{\frac {(m-t)k+r}{n}}.
\eqno (2.3)
$$

By $\alpha\in [0, 1]$ and the condition of $r$ in Grand Furuta
inequality , we can take $r\geqslant t$ if $(1-t)n\geqslant (m-t)k$,
or $r\geqslant max \{{\frac {(m-t)k-(1-t)n}{n-1}}, t\}$ if
$(m-t)k\geqslant (1-t)n$ with $n\geqslant 2$.

Take $Y(x)=(A^{\frac r 2}(A^{-{\frac t 2}}(A+xB)^{m}A^{-{\frac t
2}})^{k}A^{\frac r 2})^{\frac 1 n}$. According to (2.3), $
Y(x)\geqslant Y(0)=A^{\frac {(m-t)k+r}{n}}$ for any $x\geqslant 0$,
then $Y'(0)\geqslant 0$. Differentiate $Y^{n}(x)=A^{\frac r
2}(A^{-{\frac t 2}}(A+xB)^{m}A^{-{\frac t 2}})^{k}A^{\frac r 2}$,
use Lemma 2.1, then take $x=0$, the following equality holds:
\begin{eqnarray*}
&\ &\left. {{\frac {d}{dx}}[Y^{n}(x)]}\right |
_{x=0}=\sum^{n}_{j=1}Y(0)^{n-j}Y'(0)Y^{j-1}\\
&=&\left. {{\frac {d}{dx}}[A^{\frac r 2}(A^{-{\frac t
2}}(A+xB)^{m}A^{-{\frac t 2}})^{k}A^{\frac r 2}]}\right | _{x=0}\\
&=&A^{\frac r 2}\{\sum_{i=1}^{k}[\left. {(A^{-{\frac t
2}}(A+xB)^{m}A^{-{\frac t 2}})^{k-i}}\right | _{x=0}]\cdot[\left.
{(A^{-{\frac t
2}}(A+xB)^{m}A^{-{\frac t 2}})'}\right| _{x=0}]\\
& &\cdot[\left. {(A^{-{\frac t 2}}(A+xB)^{m}A^{-{\frac t 2}})^{i-1}}\right|
_{x=0}]\}A^{\frac r 2}\\
&=&A^{\frac r 2}\{\sum_{i=1}^{k}[A^{(m-t)(k-i)}(A^{-{\frac t
2}}(\sum^{m}_{j=1}A^{m-j}BA^{j-1})A^{-{\frac t 2
}})A^{(m-t)(i-1)}]\}A^{\frac r 2}.
\end{eqnarray*}

Replace $Y(0)$ by $A^{\frac {(m-t)k+r}{n}}$, $Y'(0)$ by $X$, we have
\begin{eqnarray*}
&\ &\sum_{j=1}^{n}A^{{\frac {(m-t)k+r}{n}}(n-j)}XA^{{\frac
{(m-t)k+r}{n}}(j-1)}\\
&=&A^{\frac r 2}\{\sum^{k}_{i=1}A^{(m-t)(k-i)}[A^{-{\frac t
2}}(\sum^{m}_{j=1}A^{m-j}BA^{j-1})A^{-{\frac t
2}}]A^{(m-t)(i-1)}\}A^{\frac r 2}.
\end{eqnarray*}

Replace $A$ by $A^{\frac {n}{(m-t)k+r}}$ in above equality, then
(2.1) is obtained. $\quad\square$

\noindent{\bf Remark 2.1.} If take $t=0$ and $k=1$ in Theorem 2.1,
this theorem is just Theorem A, which is the main result of
\cite{Positive semidefinite solutions}.

\noindent{\bf Example 2.1.} We use the same example as \cite{Positive semidefinite solutions}: For two $l\times l$ matrices $A$ and $B$, take $A=diag(a_{1}, a_{2}, \ldots ,a_{2})$, all entries of $B$ are 1. If $m$, $n$, $k$ are
positive integers, $t\in[0, 1]$, there exists positive semidefinite
matrix $X$ which satisfies:
\begin{eqnarray*}
&\ &\sum_{j=1}^{n}A^{{\frac {(m-t)k+r}{n}}(n-j)}XA^{{\frac
{(m-t)k+r}{n}}(j-1)}\\
&=&A^{\frac r 2}\{\sum^{k}_{i=1}A^{(m-t)(k-i)}[A^{-{\frac t
2}}(\sum^{m}_{j=1}A^{m-j}BA^{j-1})A^{-{\frac t
2}}]A^{(m-t)(i-1)}\}A^{\frac r 2}
\end{eqnarray*}
for $r$ such that $ \begin{cases}
        r\geqslant t ,  & \text{if }\   (1-t)n\geqslant (m-t)k \  ;\\
        r\geqslant max \{{\frac {(m-t)k-(1-t)n}{n-1}} , t\} , & \text{if }
        \   (m-t)k\geqslant (1-t)n \  \text{with } n\geqslant 2 \ .
\end{cases}$
\\

It is easy to calculate the expression of $X$:
$$
X=\Bigg(\frac {a_{p}^{\frac {r-t}{2}}a_{q}^{\frac {r-t}{2}}\Big( \Sigma^{k}_{i=1}a_{p}^{(m-t)(k-i)}a_{q}^{(m-t)(i-1)}\Big)\Big(\Sigma^{m}_{j=1}a_{p}^{m-j}a_{q}^{j-1}\Big)}
{\Sigma_{j=1}^{n}a_{p}^{\frac {((m-t)k+r)(n-j)}{n}}a_{q}^{\frac {((m-t)k+r)(j-1)}{n}}}\Bigg)_{p, q=1, 2, \ldots , l}.
\eqno (2.4)
$$

\noindent{\bf Remark 2.2.} The condition of $r$ in Theorem 2.1 is necessary.
If the condition cannot be fulfilled, the solution of the equation may be not
 positive semidefinite.

  For example, take
 $$
 A=\begin{pmatrix} 1 & 0 \\ 0 & 2 \end{pmatrix},
  B=\begin{pmatrix} 1 & 1 \\1 & 1 \end{pmatrix},
  $$
  and $m=2$, $n=2$, $k=2$, $t={\frac 1 2}$ in
  Example 2.1. If we put $r ={\frac 1 2}$, then $r\ngeqslant max \{{\frac {(m-t)k-(1-t)n}{n-1}} , t\}$.
  By (2.4), the solution of the following matrix equation
\begin{eqnarray*}
A^{\frac 7 4}X+XA^{\frac 7 4}&=&A^{\frac 1 4}\Big(A^{\frac 3 2}\big(A^{\frac 3 4}BA^{-{\frac 1 4}}+A^{-{\frac 1 4}}BA^{\frac 3 4}\big)+\big(A^{\frac 3 4}BA^{-{\frac 1 4}}+A^{-{\frac 1 4}}BA^{\frac 3 4}\big)A^{\frac 3 2}\Big)A^{\frac 1 4}\\
&=&A^{\frac 5 2 }B+A^{\frac 3 2}BA+ABA^{\frac 3 2}+BA^{\frac  5 2}\\
&=&\begin{pmatrix} 4 & 3+6\times 2^{\frac 1 2}  \\ 3+6\times 2^{\frac 1 2}  & 16\times 2^{\frac 1 2} \end{pmatrix}
\end{eqnarray*}
is
 $$
 X=\begin{pmatrix} 2 & {\frac {3+6\times 2^{\frac 1 2 }}{1+2\times 2^{\frac 3 4}}} \\ {\frac {3+6\times 2^{\frac 1 2 }}{1+2\times 2^{\frac 3 4}}} &2\times 2^{\frac 3 4}\end{pmatrix}.
 $$
 However, $X$ is not a positive semidefinite matrix because its eigenvalues are $ \{ 5.4007\ldots, -0.0372\ldots\}$.

\noindent{\bf Remark 2.3.} In \cite{Bhatia}, the authors showed that if $A$ and $Y$ are positive semidefinite matrices in matrix equation $A^{n-1}X+A^{n-2}XA+ \cdots +AXA^{n-2}+XA^{n-1}=Y$, then so is $X$. By Theorem 2.1, in some special cases, if $Y$ can be expressed as the right hand of (2.1), though it is not positive semidefinite, then there still exists positive semidefinite solution satisfies   $A^{n-1}X+A^{n-2}XA+ \cdots +AXA^{n-2}+XA^{n-1}=Y$.

For example, take
 $$
A=\begin{pmatrix} 1 & 0 \\ 0 & 2\times2^{\frac 1 3} \end{pmatrix},  Y=\begin{pmatrix} 4 & 3\times2^{\frac 1 4}+6\times2^{\frac 3 4} \\ 3\times2^{\frac 1 4}+6\times2^{\frac 3 4}  & 32 \end{pmatrix}.
$$
Although $Y$ is not a positive semidefinite matrix (because its eigenvalues are $\{37.5589\ldots, -1.5589\ldots\}$) , by simple calculation, the solution of the following matrix equation
  $$
  A^{2}X+AXA+XA^{2}=Y
 $$
is
   $$
   X=\begin{pmatrix} {\frac 4 3} & {\frac {3\times2^{\frac 1 4}+6\times2^{\frac  3 4}}{1+2\times2^{\frac 1 3}+4\times2^{\frac 2 3}}} \\ {\frac {3\times2^{\frac 1 4}+6\times2^{\frac  3 4}}{1+2\times2^{\frac 1 3}+4\times2^{\frac 2 3}}}  & {\frac {4\times2^{\frac 1 3}}{3}} \end{pmatrix},
  $$
which is still a definite matrix whose  eigenvalues are $\{2.9013\ldots, 0.1119\ldots\}$. The critical reason is that $Y$ can be expressed as follows,
$$
Y=A^{\frac 3 8}  \{\sum^{2}_{i=1}A^{{\frac 9 8}(2-i)} [A^{-{\frac {3} {16}}}(\sum^{2}_{j=1}A^{{\frac 3 4}(2-j)}BA^{{\frac 3 4}(j-1)})A^{-{\frac {3} {16}}}]    A^{{\frac 9 8}(i-1)}     \}      A^{\frac 3 8},
$$
which is the right hand of (2.1) under the condition of $m=2$, $n=3$, $k=2$, $t={\frac 1 2}$, $r=1$.

\noindent{\bf Remark 2.4.} The following question remains open: How to investigate the properties of the solution of operator equation $X^{n-1}A+X^{n-2}AX+ \cdots +XAX^{n-2}+AX^{n-1}=B$?

\begin{center}

\end{center}
\end{document}